\documentclass[12pt,reqno]{amsart}
\usepackage{mathrsfs}
\usepackage{amsfonts}
\usepackage{amsfonts}
\usepackage{amsfonts}
\usepackage{amsfonts}
\usepackage{amsfonts}
\usepackage{amsfonts}
\usepackage{amsfonts}
\usepackage{amsfonts}
\usepackage{amsfonts}
\usepackage{amsfonts}
\usepackage{amsfonts}
\usepackage{amsfonts}
\usepackage{amsfonts}
\usepackage{amsfonts}
\usepackage{amsfonts}
\usepackage{amsfonts}
\usepackage{amsfonts}
\usepackage{amsfonts}
\usepackage{a4wide}
\let\Section=\section
\def\section{\setcounter{equation}{0}\Section}

  \sf
\begin{document}
\title{}
\begin{center}
 {\Large\bf  Riemann-Stieltjes operators and multipliers }
\end{center}
\begin{center}
 {\Large\bf  on $Q_p$ spaces in the unit ball of ${\bf
C}^{n*}$ }
\end{center}
\vspace{0.5cm}
\author{}
 \begin{center}
{\small
 Ru  Peng$^{1,2}$\ \ \ \ \ \ \ \ Caiheng  Ouyang$^{2}$}
\end{center}
\vspace{0.2cm}
\begin{center}
{\footnotesize 1.\  Department of Mathematics, Wuhan University of
Technology, Wuhan 430070, China}\\
{\footnotesize 2.\  Wuhan Institute of Physics and Mathematics, The
Chinese Academy of Sciences, Wuhan 430071, China}
\end{center}

\thanks{*Supported in part by the National Natural
Science Foundation of China (No. 10971219) and the Fundamental
Research Funds for the Central Universities (No. 2010-Ia-023).}

\vspace{0.5cm}

{\footnotesize\bf Abstract.} {\footnotesize This paper is devoted to
characterizing the Riemann-Stieltjes operators and pointwise
multipliers acting on M${\rm \ddot{o}}$bius invariant spaces $Q_p$,
which unify BMOA and Bloch space in the scale of $p$. The
boundedness and compactness of these operators on $Q_p$ spaces are
determined by means of an embedding theorem, i.e. $Q_p$ spaces
boundedly embedded in the non-isotropic tent type spaces
$T_q^\infty$. }\vspace{0.1cm}

\vspace{0.1cm}


\date{}\maketitle{}

\vskip 0.1cm \arraycolsep1.5pt
\newtheorem{Lemma}{Lemma}[section]
\newtheorem{Theorem}{Theorem}[section]
\newtheorem{Definition}{Definition}[section]
\newtheorem{Proposition}{Proposition}[section]
\newtheorem{Remark}{Remark}[section]
\newtheorem{Corollary}{Corollary}[section]

\small

\begin{center}
{\small\bf\S1\ Introduction}
\end{center}
\vspace{0.3cm}\setcounter{section}{1} \setcounter{equation}{0}

Let $B=\{{z\in {\bf C}^n:| z|< 1}\}$ be the unit ball of ${\bf C}^n$
$(n>1)$, $S=\{z\in {\bf C}^n: | z|=1\}$ be its boundary. $d\upsilon$
denotes the normalized Lebesgue measure of $B$, i.e. $\upsilon
(B)=1$, and $d\sigma$ denotes the normalized rotation invariant
Lebesgue measure of $S$ satisfying $\sigma (S)=1$. Let $d\lambda
(z)={(1-{|z|}^2)}^{-n-1}d\upsilon (z)$, then $d\lambda (z)$ is
automorphism invariant, that is for any $\psi\in$ Aut$(B)$, $f\in
L^1 (B)$, we have
$$\int_{B}f(z)d\lambda(z)=\int_{B}f\circ\psi(z)d\lambda(z),$$
where Aut$(B)$ is the group of biholomorphic automorphisms of $B$.

We denote the class of all holomorphic functions in $B$ by $H(B)$ .
For $f\in H(B), z\in B$, its complex gradient and invariant gradient
are defined as
\begin{center}
$\nabla f(z)={\nabla}_zf=({\frac{\partial f}{\partial
z_1}(z)},\ldots,{\frac{\partial f}{\partial z_n}(z)})$,\ \ \ \
$\widetilde{\nabla}f(z)=\nabla(f\circ\varphi_z)(0)$,
\end{center}
where $\varphi _z$ is the M${\rm \ddot{o}}$bius transformation for
$z\in B$ , which satisfies $\varphi _z (0)=z$, $\varphi _z(z)=0$ and
$\varphi _z\circ\varphi _z=I$,
 and its radial derivative $Rf(z)=<\nabla f(z),\bar{z}>=\sum\limits_{j=1}^{n}\frac{\partial f}{\partial z_j}(z)z_j$.

We say that $f\in H(B)$ is an Bloch function if
$${\| f\|}_{\beta}=| f(0)|+\sup\limits_{z\in B} |\nabla f(z)|
(1-{|z|}^2)<\infty.$$ The collection of Bloch functions is denoted
by $\beta$. Correspondingly, $f$ is a little Bloch function, denoted
as $f\in \beta_0$ if $\lim\limits_{| z|\rightarrow 1}|\nabla f(z)|
(1-{|z|}^2)=0$.

Based on \cite{OYZ1} and referring to \cite{AXZ}, the so-called
$Q_p$ and $Q_{p,0}$ spaces in \cite{OYZ2} are defined as
\begin{equation}\label{1.1}
Q_p=\{f\in H(B):\sup\limits_{a\in
B}\int_{B}{|\widetilde{\nabla}f(z)|}^2G^p(z,a)d\lambda(z)<\infty\},
\end{equation}
it can be also written by (see Lemma 3.2 in \cite{LO3})
\begin{equation}\label{1.2}
Q_p=\{f\in H(B):\sup\limits_{a\in B}\int_{B}{|
Rf(z)|}^2{(1-{|z|}^2)}^2{(1-{|\varphi_a
(z)|}^2)}^{np}d\lambda(z)<\infty\}.
\end{equation}

$${\| f\|}_{Q_p}=| f(0)|+\sup\limits_{a\in
B}{\left(\int_{B}{|\widetilde{\nabla}f(z)|}^2G^p(z,a)d\lambda(z)\right)}^{1\over
2}$$ and
$$Q_{p,0}=\{f\in H(B):\lim\limits_{| a|\rightarrow 1}\int_{B}{|\widetilde{\nabla}f(z)|}^2G^p(z,a)d\lambda(z)=0\},$$
for $0<p<\infty$, where $G(z,a)=g(\varphi_a(z))$ and
$$g(z)={n+1\over 2n}\int_{| z|}^{1}{(1-t^2)}^{n-1}t^{-2n+1}dt.$$
About $Q_p$ and $Q_{p,0}$, the following properties are proved in
\cite{OYZ2}.

 (i) When $0<p\leq {n-1\over n}$ or $p\geq {n\over n-1}$,
$Q_p$\ ($Q_{p,0}$) are trivial, i.e. they contain only the
constant functions. When ${n-1\over n}<p<{n\over n-1}$, $Q_p$\
($Q_{p,0}$) are nontrivial.

(ii) $Q_{p_1}\subseteq Q_{p_2}$\ \ ($Q_{{p_1},0}\subseteq
Q_{{p_2},0}$) for $0<p_1\leq p_2\leq 1$.

(iii) $Q_1=BMOA$\ \ ($Q_{1,0}=VMOA$).

 (iv) $Q_p$=Bloch space\ \
($Q_{p,0}$=little Bloch space), and ${\parallel
\cdot\parallel}_{Q_p}$ is equivalent to
${\parallel\cdot\parallel}_\beta$
  for $1<p<{n\over n-1}$.

$T_g$ and $L_g$ denote the Riemann-Stieltjes operators with the
holomorphic symbol $g$ on B respectively (see \cite{Hu} and
\cite{Xi2}):
$$T_g f(z)=\int_0^1 f(tz)Rg(tz)\frac{dt}{t},\ \ \ \ \ \
\ L_g f(z)=\int_0^1 g(tz)Rf (tz)\frac{dt}{t},\ \ \ z\in B.$$

It is easy to see that the pointwise multipliers  $M_g$  are
determined by

$$M_g f(z)=g (z)f(z)=g(0)f(0)+T_g f(z)+L_g f(z),\ \ \ \ \ z\in
B.$$ Of course, in the above definition $f$ is assumed to be
holomorphic on $B$. Clearly, $T_g f=L_f g$ and the Riemann-Stieltjes
operator can be viewed as a generalization of the well known
Ces$\acute{a}$ro operator.

$T_p^{\infty}(\mu)$ denotes the non-isotropic tent type space of all
$\mu$-measurable functions $f$ on B obeying
$${\| f\|}_{{T_p^\infty}(\mu)}^2=\sup\{\delta^{-np}\int_{Q_\delta(\xi)}{|f|}^2 d \mu ; \xi\in S, \delta>0\}<\infty,$$
where $Q_\delta(\xi)=\{z\in B:|1-\langle z,\xi\rangle|<\delta\}$ for
$\xi\in S$ and $\delta >0$.

As for the Riemann-Stieltjes operators, they can be traced back to
C. H. Pommerenke's paper \cite{Po} and  A. Siskakis's paper
\cite{Si} for the Ces$\acute{a}$ro operator and the extended
Ces$\acute{a}$ro operator. Since that time, in the unit disc $D$ of
complex plan, there have been a lot of results on the
Riemann-Stieltjes operators on distinct holomorphic function spaces,
e.g. see \cite{AC}, \cite{AS}, \cite{GP}, \cite{Xi4} and the
references therein. For the case of the unit ball of ${\bf C}^n$,
recently, we can find that the research on the Riemann-Stieltjes
operators has been developing, see \cite{Hu}, \cite{LS}, \cite{Xi2}
etc.

The purpose of this paper is to study the boundedness and
compactness of the Riemann-Stieltjes operators and pointwise
multipliers on $Q_p$ spaces as an extension of J. Xiao's paper
\cite{Xi4} to the complex ball. It not only is motivated by the
importance of $Q_p$ spaces which unify BMOA and Bloch space in the
scale of $p$, but also is inspired by the good idea that a space may
be boundedly embedded in tent space as in \cite{OF4} and \cite{Xi4}.
The concept of tent space is from real harmonic analysis \cite{CMS},
however, it is indeed quick way to characterize the boundedness of
some operators acting on function spaces.

For a positive Borel measure $\mu$ on $B$, if
$${\|\mu\|}_{LCM_p}^2=\sup\{{\mu(Q_\delta(\xi))\over \delta^{np}{(\log{2\over \delta} )}^{-2}};\xi\in S, \delta>0\}<\infty,$$ we call $\mu$ a logarithmic $p$-Carleson measure;
if
$$\lim\limits_{\delta\to 0}{\mu(Q_\delta(\xi))\over \delta^{np}{(\log{2\over \delta} )}^{-2}}=0,\ \ \ {\rm for}\ \xi\in S\ {\rm uniformly},$$
we call $\mu$ a vanishing logarithmic $p$-Carleson measure. The
logarithmic Carleson measure and vanishing logarithmic Carleson
measure were introduced in \cite{ZhR} and \cite{MZ}. $\mu$ is a
usual $p$-Carleson measure if the factor ${(\log{2\over \delta}
)}^{-2}$ is deleted, and denoted by ${\|\cdot\|}_{CM_p}$ simply.

In this paper we only need to consider the case
$\frac{n-1}{n}<p<\frac{n}{n-1}$, since $Q_p$ spaces are trivial when
$0<p\leq {n-1\over n}$ or $p\geq {n\over n-1}$. The main results are
as follows.

{\bf Theorem 2.1}\ Let $\frac{n-1}{n}<p\leq q<\frac{n}{n-1}$, $\mu$
be a
 positive
Borel measure on $B$. Then the identity operator $I: Q_p\mapsto
T_q^{\infty}(\mu)$ is bounded if and only if $\mu$ is a logarithmic
$q$-Carleson measure.

{\bf Theorem 2.2}\ Let $\frac{n-1}{n}<p\leq q<\frac{n}{n-1}$, $g$ be
holomorphic on $B$, $d \mu_{q,g}(z)={|Rg
(z)|}^2{(1-{|z|}^2)}^{n(q-1)+1}dv(z)$. Then
\begin{description}
\item [(i)]\ \ $T_g: Q_p\mapsto Q_q$ is bounded if and
only if $\mu_{q,g}$ is a logarithmic $q$-Carleson measure.
\item[(ii)]\ $L_g: Q_p\mapsto Q_q$ is bounded if
and only if ${\|g\|}_{H^\infty}<\infty$.
\item[(iii)]\ $M_g: Q_p\mapsto Q_q$ is bounded if
and only if $\mu_{q,g}$ is a logarithmic $q$-Carleson measure and
${\|g\|}_{H^\infty}<\infty$.
\end{description}

{\bf Theorem 3.1}\ Let $\frac{n-1}{n}<p\leq q<\frac{n}{n-1}$, $g$ be
holomorphic on $B$, $d \mu_{q,g}(z)={|Rg
(z)|}^2{(1-{|z|}^2)}^{n(q-1)+1}dv(z)$. Then
\begin{description}
\item [(i)]\ \ $T_g: Q_p\mapsto Q_q$ is compact if and
only if $\mu_{q,g}$ is a vanishing logarithmic $q$-Carleson measure,
here the part "if" holds except for the case of ${n-1\over n}<p\leq
q<1$.
\item[(ii)]\ $L_g: Q_p\mapsto Q_q$ is compact if
and only if $g=0$.
\item[(iii)]\ $M_g: Q_p\mapsto Q_q$ is compact if
and only if $g=0$.
\end{description}

Theorem 2.1 is the base of arguments of Theorem 2.2 and Theorem 3.1.
These two theorems are extension of Theorem 1.2 of \cite{Xi4} to the
case of the unit ball of ${\bf C}^n$, especially for the operators
$T_g$, $L_g$ and $M_g$ between $Q_p$ spaces in distinct scale of $p$
and in terms of logarithmic $p$-Carleson measure defined by the
non-isotropic metric ${|1-\langle z,\xi\rangle|}^{1\over2}$ on the
ball $\bar{B}$. By Lemma 2.1 below with $s=nq$, it is easy to see
that ${\|\mu_{q,g}\|}_{LCM_q}<\infty$ with $d \mu_{q,g}(z)={|Rg
(z)|}^2{(1-{|z|}^2)}^{n(q-1)+1}dv(z)$ is equivalent to
\begin{equation}\label{1.3}
\sup\limits_{w\in B}\{\log^2{2\over {1-{|w|}^2}}\int_{B}{|
Rg(z)|}^2{(1-{|z|}^2)}^2{(1-{|\varphi_w
(z)|}^2)}^{nq}d\lambda(z)\}<\infty.
\end{equation}
Thus (iii) of Theorem 2.2 is not only an extension of Theorem 1 of
\cite{PP} to the unit ball of ${\bf C}^n$, but also an extension to
the whole range  $\frac{n-1}{n}<p\leq q<\frac{n}{n-1}$ (the results
 in \cite{PP} is only for $0<p=q<1$). Recalling another expression (1.2) of definition of
$Q_q$ spaces, the class of all symbol functions $g$ satisfying (1.3)
would be smaller than $Q_q$, which we might call a logarithmic type
$Q_q$ spaces, denoted as $\log Q_q$. In other words, the necessary
and sufficient condition ${\|\mu_{q,g}\|}_{LCM_q}<\infty$ in Theorem
2.2 may be alternatively changed into $g\in \log Q_q$, which seems
to be more convenient for verifying the boundedness of the operators
$T_g$ and $M_g$.

Among the above theorems, some new and special techniques will be
adapted to overcome the difficulty causing by the differences of one
and several complex variables or target spaces. The embedding result
for the pointwise multipliers on $Q_p$ spaces will prompt us to
solve a corona type problem for $Q_p$ spaces in the future.

Throughout this paper, $C$, $M$ denote positive constants which are
not necessarily the same at each appearance. The expression
$A\approx B$ means that there exists a positive $C$ such that
$C^{-1}B\leq A\leq CB$.
  \vspace{0.5cm}
\begin{center}
{\small\bf\S2\ \  Boundedness }
\end{center}
\vspace{0.3cm} \setcounter{section}{2} \setcounter{equation}{0}

The following lemma is a version of Lemma 3.2 of \cite{OF4} with
$q=2$, $N=s$ and replacing $n$ by $np$. We omit its proof.

\vspace{0.1cm}

 {\bf Lemma 2.1}\ \ Let $0<p<\infty$, $\mu$ be a
positive Borel measure. Then the following statements are equivalent
:

(i)\  The measure $\mu$ satisfies
$$\sup\{\mu({Q_\delta(\xi)});\xi\in
S\}\leq C\frac{\delta^{np}}{\log^2{2\over \delta}}.$$

(ii)\ For every $s>0$,
\begin{equation}\label{2.1}
\sup\{\log^2{2\over
{1-{|w|}^2}}\int_B\frac{{(1-{|w|}^2)}^s}{{|1-\langle
z,w\rangle|}^{np+s}}d\mu(z);w\in B\}<\infty.
\end{equation}

(iii)\ For some $s>0$,
$$\sup\{\log^2{2\over
{1-{|w|}^2}}\int_B\frac{{(1-{|w|}^2)}^s}{{|1-\langle
z,w\rangle|}^{np+s}}d\mu(z);w\in B\}<\infty.$$

\vspace{0.1cm}

 {\bf Lemma 2.2}\ \ Let  $n\geq 2$,
$\frac{n-1}{n}<p\leq q<\frac{n}{n-1}$, $\mu$ is a logarithmic
$q$-Carleson measure. Then, for $s>n(q-p)+1$,
$$\int_{Q_\delta(\xi)}{\left(\int_B\frac{{(1-{|w|}^2)}^s| g(w)|}{{|1-\langle
z,w\rangle|}^{n+1+s}}dv(w)\right)}^2 d\mu(z)\leq
C\delta^{nq}{\|\mu\|}_{LCM_q}^2{\left\|{|
g(z)|}^2{{(1-{|z|}^2)}^{n(p-1)-1}dv(z)}\right\|}_{CM_p}^2.$$

{\bf Proof} \ \ Let $\xi\in S$, $0<\delta\leq 2$ and
$$I_{\xi, \delta}={\left(\int_{Q_\delta(\xi)}{\left(\int_B\frac{{(1-{|w|}^2)}^s| g(w)|}{{|1-\langle
z,w\rangle|}^{n+1+s}}dv(w)\right)}^2 d\mu(z)\right)}^{1\over 2}.$$

Fix ${Q_\delta(\xi)}$, let ${\|\cdot\|}_{{Q_\delta(\xi)}}$ denote
the usual norm on $L^2({Q_\delta(\xi)},d\mu)$. By duality,
$$I_{\xi, \delta}=\sup\limits_{{\|\psi\|}_{Q_\delta(\xi)}=1}\left\{\int_{Q_\delta(\xi)}\int_B\frac{{(1-{|w|}^2)}^s| g(w)|}{{|1-\langle
z,w\rangle|}^{n+1+s}}dv(w)|\psi(z)| d\mu(z)\right\}.$$

For $j\in N$, let $A_1={Q_{4\delta}(\xi)}$ and
$A_j={Q_{4^j\delta}(\xi)}\backslash{Q_{4^{j-1}\delta}(\xi)}$, $j\geq
2$. Clearly, $B=\bigcup\limits_{j=1}^{J_\delta}A_j$, where
$J_\delta$ is the integer part of $1+\log_4{2\over \delta}$.

\begin{eqnarray*}
I_{\xi,
\delta}&\leq&\sup\limits_{{\|\psi\|}_{Q_\delta(\xi)}=1}\{\int_{Q_\delta(\xi)}\int_{Q_{4\delta}(\xi)}\frac{{(1-{|w|}^2)}^s|
g(w)|}{{|1-\langle z,w\rangle|}^{n+1+s}}dv(w)|\psi(z)| d\mu(z)\\
&+&\sum\limits_{j=2}^{J_\delta}\int_{Q_\delta(\xi)}\int_{A_j}\frac{{(1-{|w|}^2)}^s|
g(w)|}{{|1-\langle z,w\rangle|}^{n+1+s}}dv(w)|\psi(z)|
d\mu(z)\}\\
&=&\sup\limits_{{\|\psi\|}_{Q_\delta(\xi)}=1}\{I_{\xi,
\delta}^{(1)}+I_{\xi, \delta}^{(2)}\}.
\end{eqnarray*}

At first, to estimate $I_{\xi, \delta}^{(1)}$. By H${\rm
\ddot{o}}$lder's inequality and Fubini's theorem, we have
\begin{eqnarray*}
I_{\xi,
\delta}^{(1)}&\leq&{\left(\int_{Q_\delta(\xi)}\int_{Q_{4\delta}(\xi)}\frac{{|
g(w)|}^2{(1-{|w|}^2)}^s\log^2{2\over {1-{|w|}^2}}}{{|1-\langle
z,w\rangle|}^{n+1+s}}dv(w)
d\mu(z)\right)}^{1\over 2}\\
&\times&{\left(\int_{Q_\delta(\xi)}\int_{Q_{4\delta}(\xi)}\frac{{|\psi(z)|}^2{(1-{|w|}^2)}^s}{{|1-\langle
z,w\rangle|}^{n+1+s}\log^2{2\over {1-{|w|}^2}}}dv(w)
d\mu(z)\right)}^{1\over 2}\\
&=&{\left(\int_{Q_{4\delta}(\xi)}\int_{Q_{\delta}(\xi)}\frac{{|
g(w)|}^2{(1-{|w|}^2)}^s\log^2{2\over {1-{|w|}^2}}}{{|1-\langle
z,w\rangle|}^{n+1+s}} d\mu(z)dv(w)\right)}^{1\over 2}\\
&\times&{\left(\int_{Q_\delta(\xi)}\int_{Q_{4\delta}(\xi)}\frac{{(1-{|w|}^2)}^s}{{|1-\langle
z,w\rangle|}^{n+1+s}\log^2{2\over {1-{|w|}^2}}}dv(w){|\psi(z)|}^2
d\mu(z)\right)}^{1\over 2}.
\end{eqnarray*}

Similar to the proof of Lemma 3.4 in \cite{OF1}, it is clear that
the inner integral of the last line above is bounded. And note that
$1-|w|\leq |1-\langle w,\xi\rangle| < 4\delta$ for $w\in
{Q_{4\delta}(\xi)}$ and so ${(1-|w|)}^{\frac
{nq-np}{2}}<{(4\delta)}^{\frac {nq-np}{2}}$ for $p\leq q$. Therefore
\begin{eqnarray*}
I_{\xi, \delta}^{(1)}
&\leq&C{\left(\int_{Q_{4\delta}(\xi)}\int_{Q_{\delta}(\xi)}\frac{{(1-{|w|}^2)}^{s}\log^2{2\over
{1-{|w|}^2}}}{{|1-\langle z,w\rangle|}^{np+s-1}} d\mu(z){|
g(w)|}^2{(1-{|w|}^2)}^{n(p-1)-2}dv(w)\right)}^{1\over 2}\\
&\times&{\left(\int_{Q_\delta(\xi)}{|\psi(z)|}^2
d\mu(z)\right)}^{1\over 2}\\
&\leq&C\delta^{\frac{nq-np}{2}}{\left(\int_{Q_{4\delta}(\xi)}\int_{Q_{\delta}(\xi)}\frac{{(1-{|w|}^2)}^{s-nq+np-1}\log^2{2\over
{1-{|w|}^2}}}{{|1-\langle z,w\rangle|}^{nq+(s-nq+np-1)}} d\mu(z){|
g(w)|}^2{(1-{|w|}^2)}^{n(p-1)-1}dv(w)\right)}^{1\over 2}\\
&\times&{\left(\int_{Q_\delta(\xi)}{|\psi(z)|}^2
d\mu(z)\right)}^{1\over 2}.
\end{eqnarray*}

 By Lemma 2.1, we can get
\begin{eqnarray*}
\sup\limits_{{\|\psi\|}_{Q_\delta(\xi)}=1}I_{\xi,
\delta}^{(1)}&\leq&C\delta^{\frac{nq-np}{2}}{\|\mu\|}_{LCM_q}{\left(\int_{Q_{4\delta}(\xi)}{|
g(w)|}^2{(1-{|w|}^2)}^{n(p-1)-1}dv(w)\right)}^{1\over 2}\\
&\leq&C\delta^{\frac{nq}{2}}{\|\mu\|}_{LCM_q}{\left\|{|
g(w)|}^2{(1-{|w|}^2)}^{n(p-1)-1}dv(w)\right\|}_{CM_p}.\\
\end{eqnarray*}

Next to consider $I_{\xi, \delta}^{(2)}$. For $j\geq 2$, $z\in
Q_\delta(\xi)$ and $w\in A_j$, we have $${|1-\langle
w,z\rangle|}^{1\over 2}\geq{|1-\langle w,\xi\rangle|}^{1\over
2}-{|1-\langle z,\xi\rangle|}^{1\over
2}\geq{(4^{j-1}\delta)}^{1\over 2}-\delta^{1\over
2}\geq2^{j-2}\delta^{1\over 2}.$$

By these estimates, H${\rm \ddot{o}}$lder's inequality and Fubini's
theorem, we have
\begin{eqnarray*}
I_{\xi, \delta}^{(2)}
&\leq&C\sum\limits_{j=2}^{J_\delta}{(4^{j-2}\delta)}^{-n-1}\int_{Q_\delta(\xi)}\int_{A_j}|
g(w)|dv(w)|\psi(z)|
d\mu(z)\\
&\leq&C\sum\limits_{j=2}^{J_\delta}{(4^{j-2}\delta)}^{-n-1}{\left(\int_{Q_\delta(\xi)}\int_{A_j}{|
g(w)|}^2{(1-{|w|}^2)}^{n(p-1)-1}dv(w)d\mu(z)\right)}^{1\over 2}\\
&\times&{\left(\int_{Q_\delta(\xi)}\int_{A_j}{| \psi
(z)|}^2{(1-{|w|}^2)}^{1-n(p-1)}dv(w)d\mu(z)\right)}^{1\over 2}\\
&\leq&C\sum\limits_{j=2}^{J_\delta}{(4^{j-2}\delta)}^{-n-1}{\left(\int_{A_j}{|
g(w)|}^2{(1-{|w|}^2)}^{n(p-1)-1}dv(w)\right)}^{1\over 2}\\
&\times&\mu^{1\over
2}(Q_\delta(\xi))\times{\left(\int_{A_j}\left(\int_{Q_\delta(\xi)}{|
\psi
(z)|}^2d\mu(z)\right){(1-{|w|}^2)}^{1-n(p-1)}dv(w)\right)}^{1\over
2}.
\end{eqnarray*}

Therefore
\begin{eqnarray*}
\sup\limits_{{\|\psi\|}_{Q_\delta(\xi)}=1}I_{\xi,
\delta}^{(2)}&\leq&C\sum\limits_{j=2}^{J_\delta}{(4^{j-2}\delta)}^{-n-1}{(4^{j}\delta)}^{np\over
2}{\left\|{|
g(w)|}^2{{(1-{|w|}^2)}^{n(p-1)-1}dv(w)}\right\|}_{CM_p}\\
&\times&\frac{\delta^{nq\over 2}}{\log{2\over
\delta}}{\|\mu\|}_{LCM_q}{\left(\int_{A_j}{(1-{|w|}^2)}^{1-n(p-1)}dv(w)\right)}^{1\over
2}\\
&\leq&C\sum\limits_{j=2}^{J_\delta}{(4^{j-2}\delta)}^{-n-1}{(4^{j}\delta)}^{np\over
2}{\left\|{|
g(w)|}^2{{(1-{|w|}^2)}^{n(p-1)-1}dv(w)}\right\|}_{CM_p}\\
&\times&\frac{\delta^{nq\over 2}}{\log{2\over
\delta}}{\|\mu\|}_{LCM_q}{(4^j\delta)}^{2+2n-np\over 2}\\
&\leq&C\frac{J_\delta}{\log{2\over \delta}}\delta^{nq\over
2}{\|\mu\|}_{LCM_q}{\left\|{|
g(w)|}^2{{(1-{|w|}^2)}^{n(p-1)-1}dv(w)}\right\|}_{CM_p}\\
&\leq&C\delta^{nq\over 2}{\|\mu\|}_{LCM_q}{\left\|{|
g(w)|}^2{{(1-{|w|}^2)}^{n(p-1)-1}dv(w)}\right\|}_{CM_p}.
\end{eqnarray*}

Thus, we have
$$I_{\xi, \delta}\leq C\delta^{nq\over 2}{\|\mu\|}_{LCM_q}{\left\|{|
g(w)|}^2{{(1-{|w|}^2)}^{n(p-1)-1}dv(w)}\right\|}_{CM_p},$$ which
ends the proof.

{\bf Lemma 2.3}\ \ Let  $n\geq 2$, $\frac{n-1}{n}<p<\frac{n}{n-1}$.
For $w\in B$, the functions $f_w(z)=\log \frac{1}{1-\langle
z,w\rangle}$ satisfy $\sup\limits_{w\in B}{\| f_w\|}_{Q_p}<\infty$.

{\bf Proof}\ \ By Theorem 3.2 of \cite{LO3}, we have
\begin{eqnarray*}
{\| f_w\|}_{Q_p}^2&\leq&C\sup\limits_{a\in
B}\int_B{|Rf_w(z)|}^2{(1-{|z|}^2)}^{n(p-1)+1}\left({1-{|a|}^2\over
{| 1-\langle z,a\rangle|}^2}\right)^{np}dv(z)\\
&=&C\sup\limits_{a\in B}\int_B\frac{1}{{| 1-\langle
z,w\rangle|}^2}{|\langle
z,w\rangle|}^2{(1-{|z|}^2)}^{n(p-1)+1}\left({1-{|a|}^2\over
{| 1-\langle z,a\rangle|}^2}\right)^{np}dv(z)\\
&\leq&C\sup\limits_{a\in
B}{(1-{|a|}^2)}^{np}\int_B\frac{{(1-{|z|}^2)}^{n(p-1)+1}}{{|
1-\langle z,a\rangle|}^{2np}{| 1-\langle z,w\rangle|}^2}dv(z).\ \ \
\ \ \ \ \ \ \ \ \ \ \ \ \ \ \ \ \
 (2.2)
\end{eqnarray*}

Let $s=n(p-1)+1$, $r=2np$, $t=2$. It is easy to know $s>-1$,
$r,t\geq 0 $, $r+t-s>n+1$ and $t-s<n+1$. Using Lemma 2.5 of
\cite{OF2}, we have

(i)When $n\geq 3$, $r-s=np+n-1>(n-1)+(n-1)\geq n+1$,
$$(2.2)\leq
C\sup\limits_{a\in
B}{(1-{|a|}^2)}^{np}\frac{1}{{(1-{|a|}^2)}^{np-2}{| 1-\langle
a,w\rangle|}^2}\leq C\sup\limits_{a\in
B}\frac{{(1-{|a|}^2)}^2}{{(1-|a|)}^2}\leq C.$$

(ii)When $n=2$,

If $r-s<n+1$, $$(2.2)\leq C\sup\limits_{a\in
B}{(1-{|a|}^2)}^{np}\frac{1}{{| 1-\langle a,w\rangle|}^{np}}\leq
C\sup\limits_{a\in B}\frac{{(1-{|a|}^2)}^{np}}{{(1-|a|)}^{np}}\leq
C.$$

If $r-s>n+1$, $$(2.2)\leq C\sup\limits_{a\in
B}{(1-{|a|}^2)}^{np}\frac{1}{{(1-{|a|}^2)}^{np-2}{| 1-\langle
a,w\rangle|}^2}\leq C\sup\limits_{a\in
B}\frac{{(1-{|a|}^2)}^2}{{(1-|a|)}^2}\leq C.$$

If $r-s=n+1$, i.e $p=1$, by Lemma 3.1 of \cite{OF4}, we can get
$\sup\limits_{w\in B}{\| f_w\|}_{Q_1}<\infty$.

\vspace{0.5cm}

{\bf Proof of Theorem 2.1}
\vspace{0.2cm}

Suppose the identity operator $I: Q_p\mapsto T_q^{\infty}(\mu)$ is
bounded. For any $\xi\in S$ and $0<\delta<1$, we consider the
function $f_{\xi,\delta}(z)=\log\frac{2}{1-\langle
z,(1-\delta)\xi\rangle}$, by Lemma 2.6 of \cite{OF2}, we have
$$|f_{\xi,\delta}(z)|\approx\log{2\over \delta},\ \ \ \ \ \ z\in
Q_\delta(\xi),$$ and by Lemma 2.3
$$\delta^{-nq}\int_{Q_\delta(\xi)}{|f_{\xi,\delta}|}^2d\mu\leq
C{\|f_{\xi,\delta}\|}_{Q_p}^2\leq C.$$ Accordingly,
${\|\mu\|}_{LCM_q}\leq C.$

Conversely, suppose $\mu$ is a logarithmic $q$-Carleson measure. For
a holomorphic function $f$, we recall the following representation
formula
$$Rf(z)=C_\alpha\int_B Rf(w){{(1-{|w|}^2)}^\alpha\over
{(1-\langle z,w\rangle)}^{n+1+\alpha}}dv(w)$$ for $\alpha$ large
enough. Acting on the above equation by the inverse operator
$R^{-1}$,
$$f(z)=C_\alpha R^{-1}\int_B Rf(w){{(1-{|w|}^2)}^\alpha\over
{(1-\langle z,w\rangle)}^{n+1+\alpha}}dv(w),$$ and consequently, we
can get
  $$ |f(z)|\leq
C\int_B|Rf(w)|{{(1-{|w|}^2)}^\alpha\over {|1-\langle
z,w\rangle|}^{n+\alpha}}dv(w).\ \ \ \ \ \ \ \ \ \ \ \ \ \ \ \ \ \ \
\ \ \ \ \ \ \ \ \ \ \ (2.3)$$

Using (2.3) and Lemma 2.2 with $g(w)=|Rf(w)|(1-{|w|}^2)$, we have
\begin{eqnarray*}
\delta^{-nq}\int_{Q_\delta(\xi)}{|f(z)|}^2d\mu(z)&\leq&C\delta^{-nq}\int_{Q_\delta(\xi)}{\left(\int_B|
Rf(w)|\frac{{(1-{|w|}^2)}^\alpha}{{|1-\langle
z,w\rangle|}^{n+\alpha}}dv(w)\right)}^2 d\mu(z)\\
&=&C\delta^{-nq}\int_{Q_\delta(\xi)}{\left(\int_B\frac{|
Rf(w)|{(1-{|w|}^2)}{(1-{|w|}^2)}^{\alpha-1}}{{|1-\langle
z,w\rangle|}^{n+1+(\alpha-1)}}dv(w)\right)}^2 d\mu(z)\\
&\leq&C{\|\mu\|}_{LCM_q}^2{\left\|{|
Rf(z)|}^2{{(1-{|z|}^2)}^{n(p-1)+1}dv(z)}\right\|}_{CM_p}^2\\
&\leq&C{\|\mu\|}_{LCM_q}^2{\|f\|}_{Q_p}^2,
\end{eqnarray*}
the last inequality holds because the norm of $f\in Q_p$ for
$\frac{n-1}{n}<p<\frac{n}{n-1}$ is comparably dominated by the
geometric quantity
$$| f(0)|+\sup\left\{{\left(\delta^{-np}\int_{Q_\delta(\xi)}{|Rf(z)|}^2 {(1-{|z|}^2)}^{n(p-1)+1}dv(z)\right)}^{1\over 2} ; \xi\in S, \delta>0\right\}<\infty$$
by Corollary 3.2 of \cite{LO3} with $m=1$.

\vspace{0.5cm}

{\bf Proof of Theorem 2.2}
\vspace{0.2cm}

(i)\ \ Note that $R(T_g f)(z)=f(z)Rg (z)$. So, Theorem 2.1 implies
that $T_g$ maps boundedly $Q_p$ into $Q_q$ is equivalent to
${\|\mu_{q,g}\|}_{LCM_q}<\infty$.

(ii)\ \ If ${\|g\|}_{H^\infty}<\infty$, then
\begin{eqnarray*}
& &\delta^{-nq}\int_{Q_{\delta}(\xi)}{|R(L_g
f)(z)|}^2{(1-{|z|}^2)}^{n(q-1)+1}dv(z)\\
&=&\delta^{-nq}\int_{Q_{\delta}(\xi)}{|g
(z)|}^2{| Rf(z)|}^2{(1-{|z|}^2)}^{n(q-1)+1}dv(z)\\
&\leq& C{\|g\|}_{H^\infty}^2{\| f\|}_{Q_q}^2\leq
C{\|g\|}_{H^\infty}^2{\| f\|}_{Q_p}^2,
\end{eqnarray*}
this implies that ${\| L_g f\|}_{Q_q}\leq C{\|g\|}_{H^\infty}{\|
f\|}_{Q_p}.$ So, $L_g: Q_p\mapsto Q_q$ is bounded.

Conversely, suppose $L_g: Q_p\mapsto Q_q$ is bounded. We fix $\xi\in
S$ and give a point $w\in B$ near to the boundary with
$|w|>\frac{2}{3}$, there exists $0<\delta<1$ such that
$$E(w, {1\over 2})\subset Q_\delta (\xi)\ \ \ \ {\rm and}\ \ \ \ \ 1-{|w|}^2\approx\delta,$$
where $E(z,r)=\{w\in B:|\varphi_z(w)|<r\}$ denote the
pseudo-hyperbolic metric ball at $z$. Choosing
$f_w(z)=\log\frac{1}{1-\langle z,w\rangle}$. By Lemma 2.3, we know
$\sup\limits_{w\in B}{\| f_w\|}_{Q_p}\leq C$. It is well known that
$$v(E(w,{1\over 2}))\approx {(1-{|w|}^2)}^{n+1},\ \ \ 1-{|w|}^2\approx
1-{|z|}^2 \approx |1-\langle z, w\rangle|\ \ \ \ {\rm for}\ \ \ z\in
E(w,{1\over 2}).$$ Also note that for $z\in E(w,{1\over 2})$, we
have
$$1-{| \varphi_w(z)|}^2=\frac{(1-{|w|}^2)(1-{|z|}^2)}{{|1-\langle
z,w\rangle|}^2}>\frac{3}{4}.$$ Thus $$1-|\langle z,w\rangle|\leq |
1-\langle z,w\rangle|<\frac{2}{\sqrt{3}}{(1-{|w|}^2)}^{1\over
2}(1-{|z|}^2)^{1\over 2}\leq \frac{2}{\sqrt{3}}{(1-{|w|}^2)}^{1\over
2}<\frac{2}{\sqrt{3}}\cdot\frac{\sqrt{5}}{3}=\frac{2\sqrt{15}}{9},$$
this implies $|\langle z,w\rangle|>1-\frac{2\sqrt{15}}{9}$. By the
$\mathscr{M}$-subharmonicity of ${|g (w)|}^2$, we have
\begin{eqnarray*}
{|g (w)|}^2&\leq& C\frac{1}{v(E(w,{1\over
2}))}\int_{E(w,{1\over 2})}{|g (z)|}^2 dv(z)\\
 &\leq& C\frac{1}{{(1-{|w|}^2)}^{n+1}}\int_{E(w,{1\over 2})}{|g
(z)|}^2dv(z)\\
&\leq&C\delta^{-nq}\int_{E(w,{1\over 2})}\frac{{|g
(z)|}^2{(1-{|z|}^2)}^{n(q-1)+1}}{{|1-\langle z,w\rangle|}^2}dv(z)\\
&\leq&C\delta^{-nq}\int_{E(w,{1\over 2})}\frac{{|g
(z)|}^2{|\langle z,w\rangle|}^2{(1-{|z|}^2)}^{n(q-1)+1}}{{|1-\langle z,w\rangle|}^2}dv(z)\\
&\leq&C\delta^{-nq}\int_{Q_{\delta}(\xi)}{|g (z)|}^2{| R
f_w(z)|}^2{(1-{|z|}^2)}^{n(q-1)+1}dv(z)\\
&\leq& C{\|L_g (f_w)\|}_{Q_q}^2\leq C{\|L_g\|}^2
{\|f_w\|}_{Q_p}^2\leq C,
\end{eqnarray*}
and consequently, $|g (w)|\leq C$ for $|w|>\frac{2}{3}$. By maximum
modulus principle, we have $|g (w)|\leq C$ for $w\in B$. Thus $g\in
H^\infty$.

(iii)\ The "if" part follows from the corresponding ones of (i) and
(ii). We only need to see the "only if" part. Note that
$f_w(z)=\log\frac{2}{1-\langle z,w\rangle}$ belongs to $Q_p$ with
$\sup\limits_{w\in B}{\| f_w\|}_{Q_p}\leq C$ and any function $f\in
Q_p$ has the growth (see \cite{PO2})
$$|f(z)|\leq |f(0)|+C{\| f\|}_{Q_p}\log\frac{1}{1-{|z|}^2}\leq C{\| f\|}_{Q_p}\log\frac{2}{1-{|z|}^2},\ \ \ {\rm for\ every}\ z\in B.$$
So, if $M_g: Q_p\mapsto Q_q$ is bounded, then for every $w\in B$,
$$|g(z)f_w(z)|\leq C{\| M_g
f_w\|}_{Q_q}\log\frac{2}{1-{|z|}^2}\leq C{\|
M_g\|}\log\frac{2}{1-{|z|}^2},\ \ \ \ z\in B$$ and hence $|g(w)|\leq
 C\| M_g\|$(upon taking $z=w$ in the last estimate), that is,
 ${\|g\|}_{H^\infty}<\infty$, equivalently, $L_g: Q_p\mapsto Q_q$ is
 bounded by (ii). Consequently,
 $T_g f=M_g f-L_g f- f(0)g(0)$ gives the
 boundedness of $T_g: Q_p\mapsto Q_q$  and then
 ${\|\mu_{q,g}\|}_{LCM_q}<\infty$.

\vspace{0.2cm}

 {\bf Corollary 2.1}\ \ Let $1<q<\frac{n}{n-1}$, $g$ be holomorphic on B,
$d \mu_{q,g}(z)={|Rg (z)|}^2{(1-{|z|}^2)}^{n(q-1)+1}dv(z)$ and
${\|g\|}_{H^\infty}=\sup\limits_{z\in B}{| g (z)|}$. Then
\begin{description}
\item [(i)]\ \ $T_g: BMOA\mapsto \beta$ is bounded if and
only if $\mu_{q,g}$ is  logarithmic $q$-Carleson measure.
\item[(ii)]\ $L_g: BMOA\mapsto \beta$ is bounded if
and only if ${\|g\|}_{H^\infty}<\infty$.
\item[(iii)]\ $M_g: BMOA\mapsto \beta$ is bounded if
and only if $\mu_{q,g}$ is  logarithmic $q$-Carleson measure and
${\|g\|}_{H^\infty}<\infty$.
\end{description}

 \vspace{0.5cm}
\begin{center}
{\small\bf\S3\ \  Compactness  }
\end{center}
\vspace{0.3cm} \setcounter{section}{3} \setcounter{equation}{0}

\vspace{0.3cm}

Before proving the compactness of $T_g$, $L_g$ and $M_g$, we need to
give the following lemmas.

{\bf Lemma 3.1}(Lemma 3.7 of \cite{Tj})\ Let $X$, $Y$ be two Banach
spaces of analytic functions on $D$. Suppose

(1) the point evaluation functionals on $Y$ are continuous;

(2) the closed unit ball of $X$ is a compact subset of $X$ in the
topology of uniform convergence on compact sets;

(3) $T: X\mapsto Y$ is continuous when $X$ and $Y$ are given the
topology of uniform convergence on compact sets.\\
Then $T$ is a compact operator if and only if given a bounded
sequence $\{f_j\}$ in $X$ such that $f_j\rightarrow 0$ uniformly on
compact sets, then the sequence $\{Tf_j\}$ converges to zero in the
norm of $Y$.

Although this lemma is shown for the unit disc $D$ of the complex
plane, it is still valid for any complex domain, of course,
including the unit ball of ${\bf C}^n$. In this section, to prove
the compactness of the operators $T_g$ and $L_g$ from $Q_p$ to
$Q_q$, we need to verify the three assumptions of the above lemma.

At first, it is clear that the assumption (1) holds by setting
$e_z(f)=f(z): Q_q\rightarrow {\bf C}$ because any function $f\in
Q_q$ has the growth
$$|f(z)|\leq | f(0)|+C{\| f\|}_{Q_q}\log\frac{1}{1-{|z|}^2},\ \ \ z\in B.\ \ \ \ \ \ \ \ \ \ \ \ \ \ (3.1)$$

Let $\{f_j\}$ be a sequence in the closed unit ball $\mathscr{B}$ of
$Q_p$. Since the functions in $\mathscr{B}$ are bounded uniformly on
compact sets of $B$, by Montel's theorem we can pick out a
subsequence $f_{j_k}\rightarrow h$ uniformly on compact sets of $B$,
for some $h\in H(B)$. To verify the assumption (2), we show that
$h\in Q_p$. Indeed,
\begin{eqnarray*}
& &\int_{B}{|\widetilde{\nabla}h(z)|}^2G^p(z,a)d\lambda(z)\\
&=&\int_{B}\lim\limits_{k\rightarrow \infty}{|\widetilde{\nabla}f_{j_k}(z)|}^2G^p(z,a)d\lambda(z)\\
&\leq& \liminf\limits_{k\rightarrow \infty}\int_{B}{|\widetilde{\nabla}f_{j_k}(z)|}^2G^p(z,a)d\lambda(z)\\
&\leq&\liminf\limits_{k\rightarrow \infty}{\| f_{j_k}\|}_{Q_p}^2\leq
1
\end{eqnarray*}
by Fatou's lemma for every $a\in B$, so $h\in \mathscr{B}$.

The assumption (3) means that if bounded sequence $\{f_j\}$ in $Q_p$
converges uniformly to zero on compact sets of B, then $\{T_g f_j\}$
(and $\{L_g f_j\}$) converges uniformly to zero on compact sets of
B. Now we verify it. Let $f_j(z)\rightarrow 0$ uniformly on compact
sets $G$ of $B$, then $\{D^\alpha f_j\}$ converges uniformly to zero
on compact sets $K$ of $B$ and $\sup\limits_{z\in K}|D^\alpha
f_j|\leq C_\alpha \sup\limits_{z\in G}| f_j|$ by the well-known
Weierstrass theorem. Therefore
\begin{eqnarray*}
\lim\limits_{j\rightarrow \infty}| L_g
f_j(z)|&=&\lim\limits_{j\rightarrow \infty}\left|\int_0^1
g(tz)Rf_j (tz)\frac{dt}{t}\right|\\
&\leq&\lim\limits_{j\rightarrow \infty}\int_0^1
|g(tz)||\nabla f_j (tz)|dt\\
&=&\int_0^1\lim\limits_{j\rightarrow \infty}
|g(tz)||\nabla f_j (tz)|dt\\
&=& 0,\ \ \ \ \ \ {\rm uniformly\ on}\ K
\end{eqnarray*}
by Lebesgue's dominated convergence theorem, since the sequence
$\{f_j(tz)\}$ is bounded uniformly for $j$ and $t\in [0, 1]$ from
(3.1). Similarly, the assumption (3) can be verified for $T_g$.

Summarizing the above arguments, we can get a criterion of the
compactness of the $T_g$ and $L_g$ as follows.

 {\bf Lemma 3.2}\ For the Riemann-Stieltjes operators $T_g$
 and $L_g$ with the holomorphic symbol $g$, the
 following statements are equivalent

 (i) $T_g$ (resp. $L_g$) is a compact operator from $Q_p$
 to $Q_q$.

 (ii) For every bounded sequence $\{f_j\}$ in $Q_p$ such that
 $f_j\rightarrow 0$ uniformly on compact sets of B, then the sequence
 $\{T_g f_j\}$ (and $\{L_g f_j\}$) converges to zero in
 the norm of $Q_q$.

\vspace{0.2cm}

For $\xi\in S$ and $\delta>0$, set
$$Q'_\delta(\xi)=\{\eta\in S: | 1-\langle \eta,\xi\rangle|<\delta\}.$$
$Q'_\delta$ is a ball of radius $\delta^{1\over 2}$ on $S$ in the
non-isotropic metric ${| 1-\langle \eta,\xi\rangle|}^{1\over 2}$ and
note that $Q'_\delta=S$ when $\delta>2$. We have the following
covering lemma, which is a version of nonisotropic ball of Lemma
2.22 in \cite{Zh2}, and will be used in the proof of Theorem 3.1 and
elsewhere later on.

{\bf Lemma 3.3}\ \ Given any natural number $m$, there exists a
natural number $N$ such that every nonisotropic ball of "radius"
$\delta\leq 2$, can be covered by $N$ nonisotropic balls of "radius"
$\delta/m$.

{\bf Proof}\ \ The first half of the proof is the same process as
that of Lemma 2.22 of \cite{Zh2} for Bergman metric ball. We can get
a covering $\{Q'_{\delta/m}(\xi'_k)\}$ of $Q'_\delta (\xi)$ with
$|1-\langle\xi'_i,\xi'_j\rangle|\geq\delta/2m$ for $i\neq j$ where
each $\xi'_k\in Q'_\delta(\xi)$. We omit its details.

Since the sets $\{Q'_{\delta/4m}(\xi'_k)\}$ are disjoint and
contained in $Q'_{\delta+(\delta/4m)}(\xi)$, we can get
$$\sum\limits_k \sigma(Q'_{\delta/4m}(\xi'_k))
\leq\sigma(Q'_{\delta+(\delta/4m)}(\xi)).$$
Moreover, there is a
positive constant $C$, independent of $\delta$ but dependent on $m$
such that
$$\sigma(Q'_{\delta+(\delta/4m)}(\xi))
\leq C\sigma(Q'_{\delta/4m}(\xi'_k))$$ for each $k$. This inequality
is true because Proposition 5.1.4 of \cite{Ru} implies that
$$C=\sup\limits_{Q'}\frac{\sigma(Q'_{\delta+(\delta/4m)})}{\sigma(Q'_{\delta/4m})}\leq{(4m+1)}^n A_0/2^{-n}.$$
Thus we see that $k\leq C$ and so the natural number $N=[C]+1$ as
desired.

 \vspace{0.2cm}

In the proof of Theorem 3.1 (ii), we need to use the lemma below,
which is also of independent interest.

{\bf Lemma 3.4}\ \ Let $f(z)$ be bounded holomorphic function on
$B$, ${\|f\|}_{H^\infty}=\sup\limits_{z\in B}|f(z)|$. Then
$$|f(z_1)-f(z_2)|\leq 2{\|f\|}_{H^\infty}|\varphi_{z_1}(z_2)|$$
holds for $z_1, z_2\in B$, where $\varphi_z$ is the M${\rm
\ddot{o}}$bius transformation of $B$.

{\bf Proof}\ \ The conclusion for the unit disc $D$ was pointed out
in \cite{Xi4}. In fact, using the invariant form of Schwarz's lemma,
i.e. Schwarz-Pick lemma for $f(z)/{\|f\|}_{H^\infty}$, we have
$$\frac{\frac{1}{{\|f\|}_{H^\infty}}\left|f(z_1)-f(z_2)\right|}{\left|1-\frac{1}{{\|f\|}_{H^\infty}^2}\overline{f(z_2)}f(z_1)\right|}\leq
\left|\frac{z_1-z_2}{1-\overline{z_2}z_1}\right|=|\varphi_{z_2}(z_1)|.$$
Therefore
$$|f(z_1)-f(z_2)|\leq
\left|{\|f\|}_{H^\infty}-{\|f\|}_{H^\infty}^{-1}\overline{f(z_2)}f(z_1)\right||\varphi_{z_1}(z_2)|\leq
2{\|f\|}_{H^\infty}|\varphi_{z_1}(z_2)|.\ \ \ \ \ \ \ \ \ (3.2)$$

In the case of the unit ball of ${\bf C}^n$, $n\geq 2$, we consider
the function $f(z)/{\|f\|}_{H^\infty}$ and $z_1,z_2\in B$. Given a
unitary map $U$ so that $z_3=U^{-1}(\varphi_{z_1}(z_2))\in ({\bf
C},0,\cdots,0)\in {\bf C}^n$ and let $F(z)=f\circ \varphi_{z_1}\circ
U/{\|f\|}_{H^\infty}$ restricted to $(D,0,\cdots,0)$. Then
$|z_3|=|\varphi_{z_1}(z_2)|$ and $F(z_3)=f(z_2)/{\|f\|}_{H^\infty}$,
$F(0)=f(z_1)/{\|f\|}_{H^\infty}$. Since $F(z)$ is a function from
the unit disc to itself, applying Schwarz-Pick lemma to $F(z)$, we
get
$$\frac{\frac{1}{{\|f\|}_{H^\infty}}\left|f(z_1)-f(z_2)\right|}{\left|1-\frac{1}{{\|f\|}_{H^\infty}^2}\overline{f(z_2)}f(z_1)\right|}=\frac{|F(0)-F(z_3)|}{|1-\overline{F(z_3)}F(0)|}
\leq |z_3|=|\varphi_{z_1}(z_2)|.$$ This means that (3.2) is still
true for the unit ball of ${\bf C}^n$.

 \vspace{0.5cm}

{\bf Proof of Theorem 3.1}
 \vspace{0.2cm}

(i)\ Suppose $\mu_{q,g}$ is a vanishing logarithmic $q$-Carleson
measure. Let $\{f_j\}$ be any bounded sequence in $Q_p$ and
$f_j\rightarrow 0$ uniformly on compact sets of $B$. For the
compactness of $T_g$, it suffices to prove
$\lim\limits_{j\rightarrow\infty}{\|T_g f_j\|}_{Q_q}=0$ by Lemma
3.2.

For $r\in (0,1)$, define the cut-off measure
$d\mu_{q,g,r}(z)=\chi_{\{z\in B: |z|>r\}}d\mu_{q,g}(z)$, where
$\chi_E$ denotes the characteristic function of a set $E$ of $B$.
\begin{eqnarray*}
& &\int_{Q_{\delta}(\xi)}{|R(T_g f_j)(z)|}^2{(1-{|z|}^2)}^{n(q-1)+1}dv(z)\\
&=&\int_{Q_{\delta}(\xi)}{|f_j(z)|}^2{|Rg
(z)|}^2{(1-{|z|}^2)}^{n(q-1)+1}dv(z)\\
&=&\int_{Q_{\delta}(\xi)}{|f_j(z)|}^2d\mu_{q,g}(z)\\
&=&\int_{Q_\delta(\xi)}{|f_j|}^2\chi_{\{z\in B: |z|\leq
r\}}d\mu_{q,g}(z)+\int_{Q_\delta(\xi)}{|f_j|}^2\chi_{\{z\in B:
|z|>r\}}d\mu_{q,g}(z)\\
&\leq&\int_{Q_\delta(\xi)}{|f_j|}^2\chi_{\{z\in B: |z|\leq
r\}}d\mu_{q,g}(z)+C\delta^{nq}{\|f_j\|}_{Q_p}^2{\|\mu_{q,g,r}\|}_{LCM_q}^2.\
\ \ \ \ \ \ \ \ \ \ \ \ \ \ \ \ \ \ \ \ \ \ \ \ (3.3)
\end{eqnarray*}
The second term of the end of (3.3) follows from the proof of the
"if" part of Theorem 2.1.

We claim that ${\|\mu_{q,g,r}\|}_{LCM_q}\rightarrow 0$ when
$r\rightarrow 1$ for the cut-off measure in the case of $1\leq
q<{n\over n-1}$. In the proof of Theorem 4.1 of \cite{LO3}, we know
$Q_\delta(\xi)\subset\widehat{Q}_{4\delta}(\xi)\subset
Q_{16\delta}(\xi)$, where
$$\widehat{Q}_\delta(\xi)=\{z\in B:\frac{z}{|z|}\in
Q'_\delta(\xi),1-\delta<|z|<1\}.$$ Hence we can use $Q_\delta$ or
alternatively $\widehat{Q}_\delta$ in the definition of (vanishing)
Carleson type measure. For any $\varepsilon>0$, there is
$\delta_0>0$ such that
$$\mu_{q,g}(\widehat{Q}_\delta(\xi))<\varepsilon\delta^{nq}{\left(\log{2\over
\delta}\right)}^{-2}$$ for all $\delta\leq\delta_0$ and for $\xi\in
S$ uniformly, since $\mu_{q,g}$ is vanishing logarithmic
$q$-Carleson measure. If $\delta>\delta_0$, given a natural number
$m=[{\delta\over \delta_0}]+1$ ($<{2\delta\over \delta_0}$) so that
${\delta\over m}<\delta_0$ for all $\delta\leq 2$, $Q'_\delta$ can
be covered by $N$ balls $Q'_{\delta/m}$ on $S$ by Lemma 3.3.
Further, it follows that
$$\widehat{Q}_\delta\cap\{z\in B:|z|>r_0\}\subset\bigcup\limits_N\widehat{Q}_{\delta/m}$$ with
$r_0=1-{\delta_0\over m}$ from the definition of
$\widehat{Q}_\delta$. Therefore
\begin{eqnarray*}
\mu_{q,g,r_0}(\widehat{Q}_\delta)&\leq&\mu_{q,g,r_0}(\bigcup\limits_N\widehat{Q}_{\delta/m})\leq\mu_{q,g,r_0}(\bigcup\limits_N\widehat{Q}_{\delta_0})
\leq\sum\limits_N\mu_{q,g}(\widehat{Q}_{\delta_0})\\
&\leq& N\varepsilon\delta_0^{nq}{\left(\log{2\over
\delta_0}\right)}^{-2}<C\frac{2^n\delta^n}{\delta_0^n}\varepsilon
\delta_0^{nq}{\left(\log{2\over
\delta_0}\right)}^{-2}\\
&<&C\varepsilon \delta^n \delta_0^{nq-n}{\left(\log{2\over
\delta_0}\right)}^{-2}<C\varepsilon \delta^{nq}{\left(\log{2\over
\delta}\right)}^{-2},
\end{eqnarray*}
where we use $N\leq C m^n$, $C>1$ (see Lemma 3.3) and $nq-n\geq 0$
when $1\leq q<{n\over n-1}$. It is clear that
$\mu_{q,g,r_0}(\widehat{Q}_\delta)<\varepsilon
\delta^{nq}{\left(\log{2\over \delta}\right)}^{-2}$ holds for
$\delta\leq\delta_0$. As shown above, for any $\varepsilon>0$ we may
find $r_0=1-{\delta_0\over m}$, so that
$$\frac{\mu_{q,g,r}(\widehat{Q}_\delta)}{\delta^{nq}{\left(\log{2\over \delta}\right)}^{-2}}<C\varepsilon
\ \ \ \ \ \ \ \ \ \ \ \ \ \ \ \ \ \ \ \ \ \ \ \ \ \ \ \ \ \ \ \ \ \
\ \ \ \ \ \ \ \ \ \ \ \ \ \ \ \ \ \ \ \ \ \ \ \ \ \ \ \ \ \ \ \ \ \
\ \ \ \ \ \ \ \ \ (3.4)$$  provided $r>r_0$ and for all $\delta\leq
2$. This is as desired.

Note that the integral of the end of (3.3)
$$\int_{Q_\delta(\xi)}{|f_j|}^2\chi_{\{z\in B: |z|\leq
r\}}d\mu_{q,g}(z)\longrightarrow 0,$$ since $f_j\rightarrow 0$
uniformly on $\{z\in B: |z|\leq r\}$. Therefore,
$$\delta^{-nq}\int_{Q_{\delta}(\xi)}{|R(T_g
f_j)(z)|}^2{(1-{|z|}^2)}^{n(q-1)+1}dv(z)\longrightarrow 0$$ as
$j\rightarrow\infty$ by (3.4). Noting that $T_g f_j (0)=0$, we have
$\lim\limits_{j\rightarrow\infty}{\|T_g f_j\|}_{Q_q}=0$.

However, at present, we are not sure  the compactness of the
operator $T_g$ for the case of ${n-1\over n}<p\leq q<1$.

Conversely, suppose $T_g: Q_p\mapsto Q_q$ is compact. $\forall \xi
\in S$, $\delta_j\rightarrow 0$, we consider the functions
$$f_j(z)={\left(\log{2\over \delta_j}\right)}^{-1}{\left(\log{2\over
1-\langle z,(1-\delta_j)\xi\rangle}\right)}^2.$$ Note that
$|1-\langle z,(1-\delta_j)\xi\rangle|\geq\delta_j$. We have
\begin{eqnarray*}
{\| f_j\|}_{Q_p}^2&\approx&| f_j (0)|^2+ \sup\limits_{a\in
B}\int_B{|Rf_j(z)|}^2{(1-{|z|}^2)}^{n(p-1)+1}\left({1-{|a|}^2\over
{| 1-\langle z,a\rangle|}^2}\right)^{np}dv(z)\\
&\leq&C+C\sup\limits_{a\in B}{\left(\log{2\over
\delta_j}\right)}^{-2}\int_B{\left|\log{2\over 1-\langle
z,(1-\delta_j)\xi\rangle}\right|}^2\frac{{|\langle
z,(1-\delta_j)\xi\rangle|}^2}{{|1-\langle
z,(1-\delta_j)\xi\rangle|}^2}\\
&&\times{(1-{|z|}^2)}^{n(p-1)+1}\left({1-{|a|}^2\over
{| 1-\langle z,a\rangle|}^2}\right)^{np}dv(z)\\
&\leq&C+C\sup\limits_{a\in B}{\left(\log{2\over
\delta_j}\right)}^{-2}\int_B{\left(\log{2\over
\delta_j}\right)}^{2}\frac{{(1-{|z|}^2)}^{n(p-1)+1}}{{|1-\langle
z,(1-\delta_j)\xi\rangle|}^2}\left({1-{|a|}^2\over
{| 1-\langle z,a\rangle|}^2}\right)^{np}dv(z)\\
&=&C+C\sup\limits_{a\in
B}{(1-{|a|}^2)}^{np}\int_B\frac{{(1-{|z|}^2)}^{n(p-1)+1}}{{|
1-\langle z,a\rangle|}^{2np}{|1-\langle
z,(1-\delta_j)\xi\rangle|}^2}dv(z)
\end{eqnarray*}

Similar to the proof of (2.2), we can get ${\|f_j\|}_{Q_p}^2\leq C$
for all $j$. It is clear that $f_j\rightarrow 0$ uniformly on
compact sets of $B$ as $\delta_j\rightarrow 0$. Using Lemma 2.6 of
\cite{OF2}, we have $\left|\log{2\over 1-\langle z,
(1-\delta_j)\xi\rangle}\right|\approx \log{2\over \delta_j}$, $z\in
{Q_{\delta_j}(\xi)}$. If $T_g$ is compact, by Lemma 3.2, we know
that for any $\xi\in S$
\begin{eqnarray*}
{\mu_{q,g}(Q_{\delta_j}(\xi))\over {\delta_j^{nq}{(\log{2\over
\delta_j} )}^{-2}}} &\leq&
C\delta_j^{-nq}\int_{Q_{\delta_j}(\xi)}{|f_j|}^2d\mu_{q,g}\\
&=&C\delta_j^{-nq}\int_{Q_{\delta_j}(\xi)}{|f_j|}^2{|Rg
(z)|}^2{(1-{|z|}^2)}^{n(q-1)+1}dv(z)\\
&=& C\delta_j^{-nq}\int_{Q_{\delta_j}(\xi)}{|R(T_g
f_j)(z)|}^2{(1-{|z|}^2)}^{n(q-1)+1}dv(z)\\
&\leq& C{\|T_g f_j\|}_{Q_q}^2\rightarrow 0,\ \ \ \ as\
j\rightarrow\infty.
\end{eqnarray*}

(ii)\ It is enough to verify that if $L_g : Q_p\rightarrow Q_q$ is
compact then $g=0$. By Theorem 2.2 (ii), the compactness of $L_g$
implies $g\in H^\infty$. Now, assume $g$ is not identically equal to
0. According to the maximum principle, the boundary value function
$g|_ S$ cannot be identically the zero function. Accordingly, there
are a positive constant $\varepsilon$ and a sequence $\{w_j\}$ in B
near to the boundary with $|w_j|>\frac{2}{3}$ such that $|g(w_j)|>
\varepsilon$. By Lemma 3.4, we have
$$|g(z_1)-g(z_2)|\leq 2{\|g\|}_{H^\infty}|\varphi_{z_1}(z_2)|,\ \ \ \ \ z_1,z_2\in B.$$
This inequality implies that there is a sufficiently small number
$r>0$ such that $|g(z)|\geq\frac{\varepsilon}{2}$ for all $j$ and
$z$ obeying $|\varphi_{w_j}(z)|<r$. Note that each pseudo-hyperbolic
ball $E(w_j,r)=\{z\in B:|\varphi_{w_j}(z)|<r\}$ is contained in
$Q_{\delta_j}(\xi)$ with $1-{|w_j|}^2\approx \delta_j$. We consider
the functions
$$f_j(z)={\left(\log{2\over
{1-{|w_j|}^2}}\right)}^{-1}{\left(\log{2\over 1-\langle z,
w_j\rangle}\right)}^2.$$ Assume $|w_j|\rightarrow 1$. It is clear
that ${\|f_j\|}_{Q_p}\leq C$ and $f_j\rightarrow 0$ uniformly on
compact sets of $B$. Note that $\left|\log{2\over 1-\langle z,
w_j\rangle}\right|\approx \log{2\over \delta_j}$ for $z\in
{Q_{\delta_j}(\xi)}$ and $|\langle z,w_j\rangle|\geq C$ for $z\in
E(w_j, r)$. Thus
\begin{eqnarray*}
{\| L_g f_j\|}_{Q_q}^2&\geq&
C\delta_j^{-nq}\int_{Q_{\delta_j}(\xi)}{| R f_j(z)|}^2{| g (z)|}^2{(1-{|z|}^2)}^{n(q-1)+1}dv(z)\\
&\geq&C\delta_j^{-nq}\int_{Q_{\delta_j}(\xi)}{\left(\log\frac{2}{\delta_j}\right)}^{-2}{\left|\log{2\over
1-\langle z, w_j\rangle}\right|}^2\frac{{|\langle
z,w_j\rangle|}^2}{{|1-\langle
z,w_j\rangle|}^2}{| g (z)|}^2{(1-{|z|}^2)}^{n(q-1)+1}dv(z)\\
&\geq&C
\varepsilon^2\delta_j^{-nq}{\left(\log\frac{2}{\delta_j}\right)}^{-2}{\left(\log\frac{2}{\delta_j}\right)}^{2}\int_{|\varphi_{w_j}(z)|<r}\frac{{|\langle
z,w_j\rangle|}^2}{{|1-\langle
z,w_j\rangle|}^2}{(1-{|z|}^2)}^{n(q-1)+1}dv(z)\\
&\geq&C
\varepsilon^2\delta_j^{-nq}\int_{|\varphi_{w_j}(z)|<r}\frac{{(1-{|z|}^2)}^{n(q-1)+1}}{{|1-\langle
z,w_j\rangle|}^2}dv(z)\\
&\geq&C
\varepsilon^2\delta_j^{-nq}{(1-{|w_j|}^2)}^{nq-n-1}\int_{|\varphi_{w_j}(z)|<r}dv(z)\\
&\geq&C
\varepsilon^2\delta_j^{-nq}{(1-{|w_j|}^2)}^{nq-n-1}{(1-{|w_j|}^2)}^{n+1}\\
&\geq&C \varepsilon^2.
\end{eqnarray*}
However, the compactness of $L_g$ forces ${\|L_g
f_j\|}_{Q_q}^2\rightarrow 0$, and consequently, $\varepsilon=0$,
contradicting $\varepsilon>0$. Therefore, $g$ must be the zero
function.

(iii)\ Suppose now $M_g: Q_p\mapsto Q_q$ is compact. Then this
operator is bounded and hence ${\|g\|}_{H_\infty}<\infty$. Let
$\{w_j\}$ be a sequence in $B$ such that $| w_j|\rightarrow 1$, and
$$f_j(z)={\left(\log \frac{2}{1-{| w_j|}^2}\right)}^{-1}{\left(\log \frac{2}{1-\langle
z,w_j\rangle}\right)}^{2}.$$ Then ${\| f_j\|}_{Q_p}\leq C$ and
$f_j\rightarrow 0$ uniformly on any compact sets of $B$. So,
${\|M_g(f_j)\|}_{Q_q}\rightarrow 0$. Since
$$|g (z)f_j(z)|=|M_g(f_j)(z)|\leq
C{\|M_g(f_j)\|}_{Q_q}\log\frac{2}{1-{|z|}^2},\ \ \ z\in B,$$ we get
(by letting $z=w_j$)
$$|g(w_j)|\log\frac{2}{1-{|w_j|}^2}\leq
C{\|M_g(f_j)\|}_{Q_q}\log\frac{2}{1-{|w_j|}^2},$$ hence
$g(w_j)\rightarrow0$. Since $g$ is bounded holomorphic function on
$B$, it follows that $g=0$.

\vspace{0.2cm}

{\bf Remark }\ The compactness result corresponding to Corollary 2.1
can be obtained. We do not go into details.

\vspace{0.5cm}

{\bf Acknowledgement}\ \ The authors would like to express their
thanks to the referee for valuable comments and suggestions.

\end{document}